\documentstyle[amscd,12pt]{amsart}
\numberwithin{equation}{section}

\newcommand{\eps}{\varepsilon}

\newcommand{\hsum}{\mathop{\sum}_{h \in H}}

\newcommand{\iopl}{\mathop{\oplus}_{i\in I}}
\newcommand{\hopl}{\mathop{\oplus}_{h \in H}}
\newcommand{\omopl}{\mathop{\oplus}_{h \in \Omega}}
\def\boxtimes{\hbox{$\kern.5pt\sqcup\kern-8.0pt\sqcap$}
\kern-8.5pt\raise.6pt\hbox{$\times$}}
\newcommand{\om}{\Omega}
\newcommand{\evw}{E_{V,\Omega}}
\newcommand{\ev}{E_{V,\Omega}}
\newcommand{\xvw}{X_{V}}
\newcommand{\xv}{X_{V}}

\newcommand{\vin}{V_{\operatorname{in}(h)}}
\newcommand{\vout}{V_{\operatorname{out}(h)}}
\newcommand{\lvw}{\Lambda_{V}}

\newcommand{\lv}{\Lambda_{V}}

\newcommand{\D}{\cal D}
\newcommand{\pv}{{\cal P}_{V,\Omega}}
\newcommand{\qv}{{\cal Q}_{V,\Omega}}

\newcommand{\dimc}{\dim}

\newcommand{\Homc}{\operatorname{Hom}\,}

\newcommand{\inn}{{\operatorname{in}}}
\newcommand{\tr}{{\operatorname{tr}}}
\newcommand{\out}{{\operatorname{out}}}

\newcommand{\End}{{\operatorname{End}\,}}

\newcommand{\Coker}{{\operatorname{Coker}\,}}

\newcommand{\gge}{>\kern-3pt>}

\def\beq{\begin{eqnarray}}
\def\endeq{\end{eqnarray}}
\def\beqn{\begin{eqnarray*}}
\def\endeqn{\end{eqnarray*}}
\renewcommand\Im{\operatorname{Im}\,}
\theoremstyle{plain}
\newtheorem{thm}{Theorem}[subsection]

\newtheorem{lemma}[thm]{Lemma}
\newtheorem{cor}[thm]{Corollary}
\theoremstyle{definition}
\newtheorem{defn}[thm]{Definition}

\begin{document}
%\font\germ=eufm10
%\def\goth#1{\hbox{\germ#1}}
 
%\input lagr.tex
%\input quiver
%\end{document}
%\end
%\input intro
 
\title[ A note on quivers with symmetries]
{ A note on quivers with symmetries}
\author{Feng Xu}
\address[ ]{Department of Mathematics,
UCLA, CA 90024}
\thanks{I'd like to thank P.S.Aspinwall, M.Douglas, G.Lusztig,G.Moore
  for very informative
discussions via e-mail . The idea of using quivers with
symmetries is also known to  H.Nakajima and I'd like to thank
him for asking interesting questions. This work is partially
supported by NSF grant DMS-9500882.}
\maketitle
\begin{abstract}
We show that the bases of irreducible integrable
highest weight module of a non-symmetric  Kac-Moody
algebra, which is associated to a quiver with a nontrivial admissible
automorphism, can be naturally identified with a set of certain
invariant Langrangian irreducible subvarieties of 
certain varieties associated with the quiver defined by Nakajima.
In the case of non-symmetric affine or finite  Kac-Moody                
algebras, the bases can be  naturally identified with a set of certain
invariant Langrangian irreducible subvarieties of a  particular
deformation of singularities
of the moduli space of instantons over A-L-E spaces.
\end{abstract}
\section{Introduction}
In his remarkable paper [Na1], Nakajima studied among other things
 the geometry of the moduli space of instantons (for more precise statement,
see \S2.1.) over A-L-E spaces. He identified a basis of irreducible  integrable
highest weight module of an affine or finite symmetric  Kac-Moody
algebra with a set of certain  Langrangian irreducible subvarieties of
a particular
deformation of the moduli space of instantons over A-L-E spaces.
By this result, it seems natural (cf. [VW]) 
to  assume that the partition functions of certain twisted N=4 
supersymmetric Yang-Mills theory on  A-L-E spaces are 
essentially normailzed characters
of  irreducible  integrable
highest weight modules of an affine or finite symmetric  Kac-Moody
algebras (cf.[Kac]).  
By a well-known results of [KP], these character functions
have certain modular properties, i.e., the vector space spanned by
the characters of a fixed level admits an action of the modular
group, in accordance with S-duality of  N=4
supersymmetric Yang-Mills (SYM) theory in 4-dimension. This is the way that
the work of [Na1] is thought as evidence in favor of  S-duality. \par
[Na1] is only concerned with the symmetric Kac-Moody algebras. On the
other hand, the results of  [KP] also applies to non-symmetric affine 
 Kac-Moody algebras. It is a natural question to ask if one
can have a similar treatment of the representations of  non-symmetric 
affine  Kac-Moody algebras in terms of the geometry of the moduli space 
of instantons over  A-L-E spaces such that the  results of  [KP]
in the case of  non-symmetric affine
 Kac-Moody algebras can be thought as further evidence in favor of
S-duality. This question is the main motivation of our paper.\par
The idea is to consider quivers with symmetries, i.e., quivers with 
an admissible non-trivial automorphism (cf.\S1.4), a concept due to
G.Lusztig in [Lu]. In [Lu], Lusztig associated to each  non-symmetric
  Kac-Moody algebra a quiver with 
an admissible non-trivial automorphism which is denoted by $a$. 
Let us denote by  $u^a$ the  corresponding non-symmetric
  Kac-Moody algebra . 
$a$ induces a natural map between the quiver varieties constructed
in [Na1]. By using the results of [Lu] and [K], we will 
show that certain $a$-invariant (as a set) irreducible Lagrangian
subvarieties of  the quiver varieties are in one-to-one 
correspondence with a basis of irreducible  integrable
highest weight modules of $u^a$ (cf. Th.3.2.1). This result partially answers
the motivating question above (cf. Cor.3.2.2).\par
Let us describe in more details the content of this paper.\par
\S2 is a preliminary section. We follow closely 
[Lu] to set up  notations. In \S2.1 we define the 
quantized algebra $U$ and its simple 
integrable module $\Lambda_\lambda$ with domiant integral weight $\lambda$.
\S2.2 is a review on quivers, perverse sheaves and canonical bases.
In \S2.3 we recall the results of [K] which establishe a one-to-one
correspondence between the canonical bases and certain irreducible Lagrangian
subvarieties. In \S2.4 we intriduce the concept of admissible
automorphisms due to G.Lusztig. Lemma 2.4.1 shows that the  one-to-one
correspondence in \S2.3 is compatible with the action of the  admissible
automorphisms on the quiver varieties.\par
In \S3.1 we recall the definitions and results of certain modified
quiver varieties which are motivated by the study of
 moduli space of instantons over A-L-E
spaces due to Nakajima (cf.[Na1]). 
Lemma 3.1.1 is essentially the same as Lemma 11.5 of [Na1] . 
Lemma 3.1.2 is the most important lemma. It identifies
a basis of  simple integrable module of   $U$  
with  a 
symmetric Cartan datum 
with certain irreducible subvarieties as defined in \S3.1.
In the case of  $U$ with a finite
symmetric Cartan datum, Lemma 3.1.2 follows from 
\S11 of [Na1].

Th.3.1.3. extends the results in \S11 of
[Na1] to include all symmetric Cartan datum by using Lemma 3.1.2. \par
In \S3.2, we consider a quiver with a nontrivial
admissible automorphism $a$. Recall the associated non-symmetric 
 Kac-Moody algebra is denoted by $u^a$. Th.3.2.1, which 
identifies the basis of a simple integrable module of 
 $u^a$ with $a$-invariant
(as a set) irreducible Lagrangian subvarieties defined in \S3.1, follows
from Lemma 2.4.1, Lemma 3.1.2 and the results of [Lu]. \par
Cor.3.2.2 follows from Th.3.2.1
in the special case of affine quivers. 
The modular properties of naturally defined functions
in Cor.3.2.2. (cf. [KP]) suggests that these functions should 
be identified with the partition function of certain twisted N=4
SYM theory on A-L-E spaces together with an admissible automorphism,
in accordance with S-duality. Cor.3.2.2. thus partially answer the question
raised at the beginning of this introduction.\par
In \S3 we concludes with questions and speculations related to string-string
duality.
For a finite set $A$, we shall use $\sharp A$ to denote 
the number of elements
in the set $A$ in this paper.

\section{Preliminaries}
\subsection{ Definition of U}
A {\bf Cartan datum} is a pair $(I,\cdot)$ consisting of a finite set $I$
and a symmetric bilinear form $\nu,\nu':\rightarrow \nu \cdot \nu'$ 
on the free abelian
group ${\Bbb Z}[I]$, with values in $\Bbb Z$. It is assumed that:
(a) $i\cdot i \in $ {2,4,6,...} for any $i\in I$; \par
(b) $2 \frac {i\cdot j}{i\cdot i} \in$ {0,-1,-2...} for any
 $i\ne j$ in $I$.\par
Two Cartan data $(I,\cdot)$ and $(I,\circ)$ are said to be {\bf proportional}
if there exist integers $a,b\geq 1$ such that $ai\circ j = b i\cdot j, \forall
i,j\in I$.\par
A Cartan datum  $(I,\cdot)$ is said to be symmetric if
$i\cdot i =2, \forall i\in I$.\par
A Cartan datum  $(I,\cdot)$ is said to be simply laced if it is 
symmetric and $i\cdot j \in \{0,-1 \}, \forall i\neq j \in I$.\par
A Cartan datum  $(I,\cdot)$ is said to be irreducible if $I$ is
non-empty and for any $i\neq j$ in $I$ there exists a sequence $i=i_1,i_2,
...,i_n =j$ in $I$ such that $i_p\cdot i_{p+1} < 0$ for $p=1,2,...,n-1$.\par

A Cartan datum  $(I,.)$ is said to be of finite type if the symmetric
matrix $(i\cdot j)$ indexed by $I\times I$ is positive definite.\par
A Cartan datum  $(I,.)$ is said to be of affine type if it is 
irreducible and the symmetric
matrix $(i\cdot j)$ indexed by $I\times I$  is positive semi-definite,
but not positive definite.\par

A root datum of type $(I,.)$ consists, by definition, of \par
(a) two finitely generated free abelian groups $Y,X$ and a perfect bilinear
pairing $\langle,\rangle : Y \times X \rightarrow \Bbb Z$;\par
(b) an embedding $I\subset X (i\rightarrow i')$ and an embedding $
I\subset Y (i\rightarrow i)$ such that: \par
(c)  $\langle i,j'\rangle = 2 \frac {i\cdot j}{i\cdot i}$. \par
In particular, we have 
(d)  $\langle i,i'\rangle = 2$ for all $i$;\par
(e)  $\langle i,j'\rangle \in$ {0,-1,-2,...} for $i\ne j$. \par
Thus $(\langle i,j'\rangle )$ is a symmetrizable generalized Cartan matrix.
The imbeddings $(b)$ induce homomorphisms ${\Bbb Z}[I] \rightarrow Y,
{\Bbb Z}[I] \rightarrow X$; we shall often denote, again by $\nu$, the
image of $\nu\in {\Bbb Z}[I]$ by either of these  homomorphisms.\par

We shall denote by $u$ the universal enveloping algebra of the Kac-Moody
algebra $g$ associated with the  symmetrizable generalized Cartan matrix
above (cf.[Kac]). \par 
A root datum above is said to be $X$-regular (resp. $Y$-regular) if
the image of the embedding $I\subset X$ is linearly independent in $X$
(resp. the image of the embedding $I\subset Y$ is 
linearly independent in $Y$). In this paper we shall be concerned {\bf only} 
with a root datum
which is both  $X$-regular and  $Y$-regular. Such a  root datum
always exists by \S2.2.2 of [Lu]. We shall choose such a root datum
for each Cartan datum. \par
Let $v$ be an indeterminate. For any $i\in I$, we set $v_i = v^{i\cdot i/2}$.
% For $m\in Z$ and $s\in N$, we set
%$$
%\bmatrix m\\ s\endbmatrix _i
%= \prod _{t=1}^s \frac {v_i^{m-t+1}-v_i^{{-(m-t+1)}}}{v_i^t - v_i ^{-t}}
%$$
Assume that a Cartan datum $(I,\cdot)$ is given.

Let   $(Y,X,\langle,\rangle,...)$ 
 be the chosen root datum of type  $(I,\cdot)$ as above.\par 
\begin{defn}
$U$ is the associative algebra over the field ${\Bbb Q}(v)$
of rational
functions of $v$ with generators $E_i, F_i, K_\mu, i\in I,\mu\in Y
$      
and the following defining relations:
\beqn
K_0 &=& 1, K_\mu K_{\mu'} = K_{\mu +\mu'}\\
K_\mu E_j K_\mu^{-1} &=& v^{\langle i,\mu\rangle} E_j, \ \ \ 
K_\mu F_j K_\mu^{-1} =
v^{-\langle i,\mu\rangle} F_j \\ 
 E_i F_j - F_j E_i &=& \delta_ {ij} \frac{\tilde K_i -
\tilde K_i^{-1}}{v_i - v_i^{-1}}\\
\sum_{r+r'=1-\langle i,j'\rangle} & &(-1)^r 
 (X_i^\pm)^{(r)} X_j^\pm (X_i^\pm)^{(r')}
= 0 \ \ \text{\rm if} \ \ i \neq j.
\endeqn
\end{defn}
Here we used notations 
$X_i^+= E_i, X_i^-=F_i,\tilde K_{\pm i} = K_{\pm (i\cdot i/2)i}$.
For any integer $r\geq 0$, and for $1\leq i \leq n$,
define $(E_i)^{(r)} = \frac {(E_i)^r}{[r]_v!}$,
 $(F_i)^{(r)} = \frac {(F_i)^r}{[r]_v!}$ where the $v$-factorial
$[r]_v!$ is defined by $[r]_v! = [r]_v.[r-1]_v.....[1]_v$. \par
It is shown in Cor.33.1.5 of [Lu] that the above definition in terms
of generators and relations
coincides with the definition in Chap.3 of [Lu].\par
Denote by $U^-$ the ${\Bbb Q}(v)$ subalgebra of $U$ generated by 
$F_i, \forall i\in I$.  \par
For $\nu\in  {\Bbb N}[I]$, we denote by $U^-_\nu$ the finite
dimensional $ {\Bbb Q}(v)$-subspace of $ U^-$ spanned by the 
monomials $F_{i_1},...,F_{i_r}$ such that for any $i\in I$, the number of 
occurences of $i$ in the sequence $i_1,...i_r$ is equal to $\nu_i$.
$U^-_\nu$ is defined in a similar way.\par
 
Define a category $C$ as follows. An object of 
$C$ is a $U $-module $M$ with a given direct sum decomposition 
$M=\oplus _{\lambda \in X}M^\lambda $ as a ${\Bbb Q}(v)$ vector space such 
that $K_i m = v^{\langle \mu , \lambda \rangle }m$ for any 
$i , \lambda $ and $m\in M^\lambda $.
 
 An object $m\in C$ is said to be 
integrable if for any $m\in M$ and $i$, there exists $n_0\geq 1$ such that
$E_i^{(n)}m=F_i^{(n)}m=0$ for all $n\geq n_0$.  
Denote by $U^- $ (resp. $_AU^- $) the subalgebra of $U $ 
(resp. $_AU$) generated by
$F_i$, $K_i^{\pm 1}$ (resp. $F_i^{(r)}, K_i^{\pm 1}$).
  The Verma module $M_\lambda $ 
is defined to be $U /J$ where
$J = \sum _i UE_i + \sum_\mu U (k_\mu - 
v^{\langle\mu , \lambda \rangle })$ 
is the left ideal of $U$.  We define 
$X^+ = \{\lambda \in X \mid \langle i, \lambda \rangle \in N ,\ \forall i\}$.
The elements of $X^+ $ are called dominant integrable weights.\par
Let $I_\lambda $ be the left ideal in $U^-$ generated by the elements 
$F_i ^{\langle i, \lambda \rangle + 1}$.
Then $\Lambda _\lambda = U^-/I_\lambda 
$ is an integrable simple 
module of $U$ 
(see Prop. 3.5.6 of [Lu]).\par
 Let $\eta_\lambda$
be the unique (up to a constant) highest vector in  $\Lambda_\lambda$.
  For $\nu \in {\Bbb N}[I]$, we define
$\Lambda_\lambda(\nu)$ to be the subspace of  $\Lambda_\lambda$ spanned by
$U_\nu^-.\eta_\lambda$.  It is clear that $\Lambda_\lambda(\nu)$ is 
the weight space of  $\Lambda_\lambda$ with weight $\lambda - \nu$.
Remember that we have used the same $\nu$ to denote its image
in $X^+$ as explained above (after the definition of root datum). \par

We can define the simple integrable module $L(\lambda)$ of $u$ with
dominant integrable weight $\lambda$ similarly as above. Let $\eta_\lambda'$
be the highest vector in  $L(\lambda)$.  For $\nu \in {\Bbb N}[I]$, we define
$L(\nu, \lambda)$ to be the subspace of  $L(\lambda)$ spanned by
$u_\nu^-.\eta_\lambda'$.\par

\subsection{ Quivers, perverse sheaves and canonical bases}
By definition, a finite {\bf graph} is a pair consisting of two finite sets 
$I$ (vertices) and $H$ (edges) and a map which to each $h\in H$ associates
a two-element subset $[h]$ of $I$.\par
Suppose a finite graph is given. In this graph, two different vertices
may be joined by several edges, but any vertex is not joined with
itself by any edges. Suppose two maps, $H\rightarrow I$ denoted by
$h\rightarrow \out(h)$,  $H\rightarrow I$ denoted by
$h\rightarrow \inn(h)$, and an involution $h \rightarrow \bar h$
are given. We assume that they satisfy the following conditions:
$$
in(\bar h) = \out (h), \out(\bar h) = \inn (h)
$$, and $\out(h) \neq \inn (h)$ for all $h\in H$.\par
An {\bf orientation} of the graph is a choice of a subset
$\Omega \subset H$ such that 
$$
\Omega \cup \bar \Omega = H, \Omega \cap \bar \Omega = \emptyset
$$.\par
A {\bf quiver} is a graph with an orientation. \par
To a graph $(I,H)$ we associate a 
Cartan datum as follows: for any $i,j\in I, i\neq j$, we define
$i\cdot j = \sharp\{h\in H~;~\out(h)=i,~\inn(h)=j\}$.
 Notice that such a  Cartan datum  is always
symmetric.\par

 We denote by $u$ the 
universal enveloping algebra of the  corresponding
Kac-Moody Lie algebra and $U$ the quantized algebra associated to the 
root datum.
%Let $k$ be an algebraic closure of a finite field $F_p$ with 
%$p$ elements.
%Let ${Cal V}$ be the category of finite dimensional $I$-graded
%$k$-vector spaces $V= \oplus_{i\in I} V_i$; the morphisms
%in  ${Cal V}$ are isomorphisms of vector spaces compatible
%with the grading.\par
Let $\cal V$ be the family of $I$-graded complex vector spaces 
$V=\iopl V_i$.

For each $V\in {\cal V}$, dim$V:= \sum_i$dim$V_i i \in {\Bbb N}[I]$
is called the dimension vector.  By abuse of notation, if  $V\in {\cal V}$
and we write $V\in  {\Bbb N}[I]$, it is assumed that we have denoted
 dim$V$ by $V$.\par
For $\nu=\sum_i \nu_i i \in {\Bbb N}[I]$, we shall sometimes identify
$\nu$ as a domiant integral weight, denoted again by $\nu$, such that
$\langle i,\nu \rangle = \nu_i$. In this way,  dim$V$ can
be identified as a  domiant integral weight.\par
For $\nu \in {\Bbb N}[I]$,
let ${\cal V}_{\nu}$ be the family of 
$I$-graded complex vector spaces $V$ 
with $\dim V=\nu$.
%Then each member of $\cal V$ belongs to ${\cal V}_{\nu}$ for
%a unique $\nu$. 
 
Let us define the complex vector spaces $\ev$ and $\xv$ by
\begin{eqnarray*}
\ev&=&\omopl \Homc( V_{\out(h)},
 V_{\inn(h)}),\\
 \xv& =&\hopl \Homc(\vout ,\vin).
\end{eqnarray*}
In the sequel, a point of $E_{V,\Omega}$ or $X_{V}$ will be denoted
as $B=(B_h)$. Here $B_h$ is in $\Homc ( V_{\out(h)},V_{\inn(h)})$.
 
We define the symplectic form $\omega $ on $X_{V}$ by
$$ \omega (B,B')=
 \hsum \eps (h)\tr(B_{\bar{h}}B'_{h}),$$
where $\eps (h)=1$ if $h\in \om$, $\eps (h)=-1$ if $h\in 
\bar{\om}$. We sometimes identify $\xvw$ and the cotangent bundle of $\evw$ 
via $\omega$.
 
The group $G_V=\prod_{i\in I}GL(V_i)$ acts on $\evw$ and $\xvw$ by
\begin{eqnarray*}
G_V\ni g=(g_i)\,&:&\,
(B_{h})\mapsto(g_{\inn(h)}B_{h}g_{\out(h)}^{-1}),
\end{eqnarray*}
where $g_i\in GL(V_i)$ for each $i\in I$.
 
The Lie algebra of $G_V$ is ${\frak g}_V=\bigoplus_{i\in I}\End(V_i)$. 
We denote an element of ${\frak g}_V$ by 
$A=(A_i)_{i\in I}$ with $A_i\in\End(V_i)$. 
The infinitesimal action of $A\in {\frak g}_V$ on $\xvw$ 
at $B\in \xvw$ is given by
$[A,B]$.
%where
%\begin{eqnarray*}
%B'_{h}&=&A_{\out(h)}B_{h}-B_{h}A_{\inn(h)},\\
%\end{eqnarray*}
Let  $\mu :\xvw \to {\frak g}_V$ be the moment map
associated with the $G_V$-action on the symplectic vector space
$\xvw$. Its $i$-th component $\mu_i:\xvw \to \End(V_i)$
is given by
\[\mu_i (B)=\sideset{}{}\sum
\begin{Sb}
h\in H \\ i=\out(h) 
\end{Sb} \eps(h)B_{\bar{h}}B_{h}.\]
 
For a non-negative integer $n$, we set 
$${\frak S}_n=\{ \sigma =(h_1,h_2,
\cdots,h_n)\,;\,h_i\in H,
\inn(h_1)=\out(h_2),\cdots,\inn(h_{n-1})=\out(h_n)\}\,,$$
 and set
${\frak S}=\bigcup_{n\ge 0}{\frak S}_n$.
For $\sigma =(h_1,h_2,\cdots,h_n)$, we
set $\out(\sigma)=\out(h_1),\inn(\sigma)=\inn(h_n)$.
For $B\in X_V$ we set
$B_{\sigma}=B_{h_n}\cdots B_{h_1}
:V_{\out(h_1)}\to V_{\inn(h_n)}$.
If $n=0$, we understand that
${\frak S}_n=\{1\}$ and $B_1$ is the identity.
An element $B$ of $X_V$ is called nilpotent if there exists a
positive integer $n$ such that $B_{\sigma}=0$
for any $\sigma\in{\frak S}_n$.
 
\begin{defn}
We set
$$X_0{}_V=\{B\in \xvw\,;\,\mu(B)=0\}$$
and
$$\lvw =\{ B\in \xvw\,;\,\mu(B)=0\text{ and } B\text{
is nilpotent}\}.$$
%Especially we put $\lv =\Lambda_{V,0}$.
\end{defn}
It is clear that $\lvw$ is a $G_V$-stable 
closed subvariety of $\xvw$.
It is known that $\lvw$ is a Lagrangian variety (cf. Remark (1)
of 5.11 in [Na1]). \par
Let us recall the results of Lusztig on canonical bases. We write
$D(X)$ for the bounded derived category of complexes of  sheaves
of ${\Bbb C}$-vector spaces on the associated complex variety $X$ over 
 ${\Bbb C}$. Objects of $D(X)$ are referred to as complexes. We shall
use the notations of [BBD]; in particular, $[d]$ denotes a shift
by $d$ degrees, and for a morphism $f$ of algebraic varieties,
$f^*$ denotes the inverse image functor, $f_!$ denotes direct 
image with compact support, etc. \par
We fix an orientation $\Omega$ of quiver. Let $\nu \in   {\Bbb N}[I]$
and let $S_\nu$ be the set of all pairs $(i,a)$ where 
$i=(i_1,...,i_m)$ is a sequence of elements of $I$ and
$a=(a_1,...,a_m)$ is a sequence of non-negative integers such that
$\nu_{i} = \sum_{l, i_l=i} a_l$. Let  $(i,a) \in S_\nu.$ A flag of 
type $(i,a)$ is, by definition, a sequence 
$\phi = (V=V^0 \supset V^1 \supset ... \supset V^m = 0)$ of $I$-graded
subspace of $V$ such that, for any $l=1,2,...,m,$ the 
$I$-graded vector space $V^{l-1}/V^l$ is zero in degreed $\neq i_l$ and
has dimension $a_l$ in degree $i_l$. We define a variety 
 $\tilde{\cal F}_{\text{i},
\text{a}}$ of all pairs $(B,\phi)$ such that $B\in \ev$ and 
$\phi$ is a $B$-stable flag of type $(\text{i},
\text{a})$. The group $G_V$ acts on 
$\tilde{\cal F}_{\text{i},\text{a}}$ in natural way.
We denote by $\pi_{\text{i},\text{a}}:
\tilde{\cal F}_{\text{i},\text{a}}\to \ev$ the 
natural projection. We note that $\pi_{\text{i},\text{a}}$
is a $G_V$-equivariant proper morphism. We set $L_{\text{i},
\text{a};\om}=(\pi_{\text{i},\text{a}})_!
(1)\in \D(\ev)$.
Here $1\in \D(\tilde{\cal F}_{\text{i},\text{a}})$ 
is the constant sheaf on 
$\tilde{\cal F}_{\text{i},\text{a}}$.
By the decomposition theorem [BBD], $L_{\text{i},\text{a}
;\om}$ is a semisimple complex.
Let $\pv$ be the set of isomorphism class of simple perverse sheaves $L$ on
$\ev$ such that $L[d]$ appears as direct summand of 
$L_{\text{i},\text{a};\om}$ for some 
$(\text{i},\text{a})\in S_{\nu}$ and
some $d\in {\Bbb Z}$. We write $\qv$ for the subcategory of $\D(\ev)$
consisting of all complexes that are isomorphic to finite direct sums of
complexes of the form $L[d]$ for various simple perverse sheaves $L\in \pv$
and various $d\in {\Bbb Z}$. Any complex in $\qv$ is semisimple and
$G_V$-equivariant.
Take $V\in {\cal V}_{\nu},V'\in {\cal V}_{\nu'}
,\bar{V}\in {\cal V}_{\bar{\nu}}$ for ${\nu}={\nu'}
+{\bar{\nu}}$ in ${\Bbb N}[I]$. We consider the diagram
\beq
&&E_{\bar V,\om}\times E_{V',\om}\mathop{\longleftarrow}^{p_1}E'
\mathop{\longrightarrow}^{p_2}E''\mathop{\longrightarrow}^{p_3}\ev\,.
\label{base}
\endeq
Here $E'$ is the variety of
$(B,\bar{\phi},\phi')$ where $B\in E_{V,\Omega}$ 
and $0\to \bar V\buildrel \bar\phi\over \longrightarrow
V\buildrel\phi'\over\longrightarrow V'\to 0$
is a B-stable exact sequence of $I$-graded vector spaces,
and $E''$ is the variety of $(B,C)$
where $B\in E_{V,\Omega}$ and $C$ is a $B$-stable $I$-graded subspace of
$V$ with $\dim C=\bar\nu$.
The morphisms $p_1$, $p_2$ and $p_3$ are defined by
$p_1(B,\bar\phi,\phi')=(B|{}_{\bar V},B|_{V'})$,
$p_2(B,\bar\phi,\phi')=(B,\Im(\bar\phi))$ and $p_3(B,C)=B$.
Note that
$p_1$ is smooth with connected fiber, $p_2$ is a principal $G_{V'}\times
G_{\bar{V}}$-bundle, and $p_3$ is proper.\par
Let $L'\in {\cal Q}_{V',\om}$ and $\bar{L}\in {\cal Q}_{\bar{V},\om}$. 
Consider the exterior tenser product $\bar L\boxtimes L'$.
Then there is $(p_2)_{\flat}p_1^*(\bar L\boxtimes L')\in\D(E'')$
such that  $(p_2)^*(p_2)_{\flat}p_1^*(\bar L\boxtimes L')
\cong p_1^*(\bar L\boxtimes L')$.
We define $L'*\bar{L}\in \qv$ by 
$(p_3)_!(p_2)_{\flat}p_1^*(\bar L\boxtimes L')[d_1-d_2]$ 
where $d_i$ is the fiber dimension of $p_i$ ($i=1,2$). \par
Let ${\cal K}_{V,\om}$ be the Grothendieck group of $\qv$.
We considered it as a ${\Bbb Z}[q,q^{-1}]$-module by 
$q(L)=L[1]$, $q^{-1}(L)=L[-1]$.
Then ${\cal K}_{\om}=\mathop{\oplus}_{\nu\in Q_-}{\cal K}_{V(\nu),\Omega}$
has a structure of an associative graded ${\Bbb Z}[q,q^{-1}]$-algebra
by the operation $*$.
We denote by $F_i\in{\cal K}_{V(i),\om}$ 
the element attached to $1\in \D(E_{V(i),\om})$ where dim$V(i)=\sum_j 
\delta_{ij} j$. The following theorem is due to Lusztig (cf. Th.13.2.11
 of [Lu]):
\begin{thm}
There is a unique ${\Bbb Q}(q)$-algebra isomorphism 
$$\Gamma_{\om}:U_q^-({\frak g})\to {\cal K}_{\om}
\mathop{\otimes}_{{\Bbb Z}[q,q^{-1}]}{\Bbb Q}(q)$$
such that $\Gamma_{\om}(f_i)=F_i$.
\end{thm}

Let us identify $L\in \pv$ with $L\otimes 1\in {\cal K}_{\om}
\mathop{\otimes}_{{\Bbb Z}[q,q^{-1}]}{\Bbb Q}(q)$. We set $\text{\bf B}
=\Gamma_{\om}^{-1}(\bigsqcup_{V\in {\cal V}}\pv)$ and call it
the canonical basis of $U_q^-({\frak g})$. 

The canonical bases  $\text{\bf B}$ have some remarkable 
properties. Given $i\in I$, we set  $\text{\bf B}_{i;\geq n}= 
\text{\bf B} \cap F_i^n U^-$ and
 $^\sigma \text{\bf B}_{i;\geq n}=
\text{\bf B} \cap  U^- F_i^n $. Let $\text{\bf B}_{i;n}
=\text{\bf B}_{i;\geq n} - \text{\bf B}_{i;\geq n+1}$
and 
 $^\sigma \text{\bf B}_{i; n}= ^\sigma \text{\bf B}_{i;\geq n} -
^\sigma \text{\bf B}_{i;\geq n+1}$. By 14.3 of [Lu] we have a partition
$\text{\bf B}_V = \sqcup_{n\geq 0} \text{\bf B}_{V;i;n}$ where
$\text{\bf B}_V := \Gamma_{\om}^{-1}(\pv)$ and 
$\text{\bf B}_{V;i;n} := \text{\bf B}_V \cap \text{\bf B}_{i;n}$.
Let $\lambda \in X^+$ be a dominant integral weight and let
$\Lambda_\lambda$ be the $U$-module as defined in \S1.1. Let $\eta_\lambda$
be the image of $1\in U^-$ under the quotient map from $U^-$ to
$\Lambda_\lambda$. Let $\text{\bf B}(\lambda) = \cap_{i\in I} 
({\cup_{n; 0 \leq n \leq \langle i, \lambda \rangle}} \  ^\sigma 
\text{\bf B}_{i,n})$. Then the map $b\rightarrow b.\eta_\lambda$ define 
a bijection from  $\text{\bf B}(\lambda)$ onto a basis of 
$\Lambda_\lambda$; moreover, if $b\in \text{\bf B} - \text{\bf B}(\lambda)$,
then $b.\eta_\lambda = 0$ (cf Th.14.4.11 of [Lu]). \par
An equivalent statement is that
$$
{\cup_{i,n; n\geq \langle i, \lambda \rangle +1}} \  ^\sigma \text{\bf B}_{i;n}
$$
is a basis of $\sum_i U^- F_i^{\langle i, \lambda \rangle +1}$ which
follows from Th.14.3.2 and \S14.4.1 of [Lu].
\subsection{ Lagrangian construction of Canonical bases}
Recall that $B(\infty)$ is the crystal base of $U^-$ as defined in \S2.3
of [K].
By [GL], $\text{\bf B}$ and $B(\infty)$ are canonically identified.
For $b\in B(\infty)$
the corresponding perverse sheaf is denoted by
$L_{b,\om}$.
Let  $B(\infty; V)$ be the set of irreducible components of
$\Lambda_V$. It is shown in [K] that the set $\sqcup_V B(\infty; V)$
has a crystal structure and are isomorphic to
 $B(\infty)$ (cf. Th.5.3.2 of [K]). We denote by 
$\Lambda_b\in\mathop{\bigsqcup}_{\nu\in Q_-}B(\infty;\nu) $ 
the corresponding element to $b\in B(\infty)$ under this isomorphism. 
Thus we have a one-to-one and surjective map $\psi: \sqcup_V B(\infty; V)
\rightarrow \text{\bf B}$ such that $\psi(\Lambda_b) = L_{b,\om}$.
We shall also use the notation $\Lambda_b(V)$ (resp.  $L_{b,\om}(V)$)
if  $\Lambda_b$ (resp.  $L_{b,\om}$) is a subset of $\Lambda_V$ (resp.
belongs to  $\text{\bf B}_V$). It is clear from the definition
that  $\psi(\Lambda_b(V)) = L_{b,\om}(V)$.
By \S2.2 we have 
dim$U^-_V = \sharp \{ L_{b,\om}(V) \}$, together with (b) of
 Th.33.1.3 of [Lu]
we  conclude that: \par 
dim$u^-_V = \sharp \{ \Lambda_b(V)|\Lambda_b(V) \in
Irr{\Lambda_V} \}$, where $Irr{\Lambda_V}$ is the set of irreducible
components of $\Lambda_V$.
\par

For $B\in X_V$, we define
$$\eps_i(B)=\dimc\Coker\Big(\mathop{\oplus}_{h;\inn(h)=i}
V(\nu)_{\out(h)}@>{(B_{h})}>>V(\nu)_i\Big).$$

Let $Y$ be a smooth algebraic variety. For any $L\in \D(Y)$,
we denote by $SS(L)$ the singular support 
(or the characteristic variety) of $L$.
It is known that $SS(L)$ is a closed Lagrangian subvariety of $T^*Y$
(See \cite{KS}).
 
We recall that $T^*\ev$ is identified with $X_V$. We record the result
 derived above and Th.6.2.2 of [K]
in the following theorem:
\begin{thm}\label{estch}
\begin{enumerate}
\item   $\dimc u^-_V = \sharp \{ \Lambda_b(V)|\Lambda_b(V) \in
Irr{\Lambda_V} \}$
\item For any $L\in \pv$ the singular support $SS(L)$ is a union of 
     irreducible components of $\lv$.
\item For any $b\in B(\infty)$ and $i\in I$, we have
\end{enumerate}
\beq
&&\Lambda_b\subset SS(L_{b,\om})\subset \Lambda_b\cup 
\mathop{\bigcup}_{\eps_i(b')>\eps_i(b)}\Lambda_{b'}.
\label{char.est}
\endeq
\end{thm}
\subsection{ Admissible automorphisms of a quiver}
An {\bf admissible automorphism} of a the graph $(I,H)$ consists,
by definition, of a permutation $a: I\rightarrow I$ and a permutation
$a:H\rightarrow H$ such that for any $h\in H$, we have $[a(h)] = a[h]$
as subsets of $I$ and such that there is no edge joining two vertices
in the same $a$-orbit. We assume that there is an integer $n\geq 1$ such
that $a^n = 1$ both on $I$ and on $H$.\par
An orientation $\Omega$ of our graph is said to be compatible with 
$a$ if, for any $h\in \Omega$, $a(\inn (h)) = \inn a(h)$ and
 $a(\out(h)) = \out a(h)$. Recall that a quiver is a graph with 
an orientation.  A quiver with a symmetry, by definition, is a
quiver together with a nontrivial  admissible automorphism which is 
 compatible with the  orientation. \par
A quiver is called affine or finite A-D-E if the underlying graph
is an  affine or finite A-D-E graph. All the  affine or finite A-D-E
quivers with symmetries have been classified in \S14.1.5 and \S14.1.6
of [Lu].\par
There is a very close connection between Cartan datum and graphs
with admissible automorphisms. Given an  admissible automorphism
$a$ of a finite graph $(I,H)$, we define $I^a$ to be the set of 
$a$-orbits on $I$. For $\tilde i,\tilde j \in I^a$, we define
$\tilde i \cdot \tilde j \in {\Bbb Z}$ as follows: if $\tilde i\neq \tilde j$, 
then
$\tilde i \cdot \tilde j $ is $-1$ times the number of edges which join
some vertex in the $a$-orbit $\tilde i$ to some vertex in the 
 $a$-orbit $\tilde j$; $\tilde i \cdot \tilde i$ is $2$ times
the number of vertices in the  $a$-orbit $\tilde i$. As shown in 
\S13.2.9 of [Lu],  $(I^a, \cdot)$ is a Cartan datum. We
shall denote  the 
 universal enveloping algebra of the corresponding Kac-Moody
algebra $g^a$ by $u^a$. Conversely,
every   Cartan datum is obtained by the construction above 
from  graphs
with admissible automorphisms (cf Prop14.1.2 of [Lu]). \par
By the results of \S14.1.5 and \S14.1.6
of [Lu], all non-symmetric affine or finite Cartan datum can
be obtained from  affine or finite A-D-E graphs with a nontrivial
 admissible automorphism by the construction above up to
proportionality. \par
Recall  ${\cal V}$ is the category of finite dimensional $I$-graded
$\Bbb C$-vector spaces. For $V\in {\cal V}$, we define $a(V)$ such that
$a(V)_i = V_{a(i)}$. For $B\in X_V$, we define $a(B) \in  X_{a(V)}$
by  setting $(a(B))_h: a(V)_i \rightarrow  a(V)_j$ to be 
$(a(B))_h = B_{a(h)}$ for any $h$ with $in(h)=j, \out(h)=i$.\par
Let 
$$
{\Bbb N}[I]^a = \{\nu \in {\Bbb N}[I]| \nu_i =  \nu_{a(i)}, \forall i\in I
\}$$.
Notice that ${\Bbb N}[I]^a$ can be identified with  ${\Bbb N}[I^a]$
in a natural way. 
If $\nu\in {\Bbb N}[I]^a$, we
shall denote, again by $\nu$,
the corresponding element in
 ${\Bbb N}[I^a]$.\par
Let $U^a$ be the quantized algebra associated with Cartan datum 
$(I^a,\cdot)$ as in \S2.1. For  $\nu\in {\Bbb N}[I]^a$ and a dominant
integral weight $\lambda$, 
recall that $(U^a)^-_\nu,(u^a)^-_\nu, L(\nu,\lambda)$ are defined as in \S2.1.  
We shall chose an orientation $\Omega$ of our graph such that 
it is compatible with $a$. Such a choice always exists by \S12.1.1
of [Lu]. \par

\begin{lemma}:
Let $\psi: \Lambda_b \rightarrow L_{b,\Omega}$ be the 
map as defined in \S2.3. Then
$$
\psi(a^{-1}  \Lambda_b(V)) = a^*  (L_{b,\Omega}(V))
$$
\end{lemma}

\begin{pf}
Suppose 
$$
 a^*  (L_{b,\Omega}(V)) = L_{a_1(b), \Omega} (a^{-1}(V)),
a^{-1}  \Lambda_b(V)= \Lambda_{a_2(b)} (a^{-1}(V))
$$. Since $a: E_{V,\Omega}\rightarrow E_{a(V),\Omega}$ is an isomorphism, both
$a_1$ and $a_2$ are one-to-one and onto maps. So:
\beqn
SS( L_{a_1(b), \Omega} (a^{-1}(V))) &=& SS( a^*  (L_{b,\Omega}(V))) \\
&=& a^{-1} (SS(L_{b,\Omega}(V))) \\
&\supset&  a^{-1} (\Lambda_b(V)) \\
&=& \Lambda_{a_2(b)} (a^{-1}(V))
\endeqn
where in the third step we have used (3) of Th.2.3.1.
So we have 
$$
SS(L_{b,\Omega}(V)) \supset \Lambda_{a_2a_1^{-1}(b)}(V)
$$. Note $\epsilon_i(b)$ is bounded above by dim$V_i$. By descending
induction on $\epsilon_i(b)$, (3) of 
Th.2.3.1 implies that $a_2a_1^{-1}(b) =b$
(cf Page 18 of [K]). So  $a_2 = a_1$ and the lemma is proved.
\end{pf}
Denote by $B(\infty; V)^a$ the set of irreducible components
of $\Lambda_V$ which is $a$-invariant as a set.  Then one can show that
$\sqcup_{V\in N[I]^a} B(\infty; V)^a$ has a crystal structure which
is isomorphic to the crystal bases $B(\infty)^a$ of $(U^a)^-$. But we
shall not use this fact in the following.
\section{ Modified varieties associated with quivers [Na1]}
\subsection{}
In this section we recall the modification of the definition of
varieties associated with quivers in \S2.2 due to Nakajima (cf. [Na1]).
The motivation of such a modification comes from the study of 
moduli space of instantons over A-L-E spaces. \par
Fix a finite graph $(I,H)$ as in \S1.3. Suppose we are given pairs of 
finite dimensional hermitian vector spaces $V_k, W_k$ for each
vertex $k$. Let us define a complex vector space $M(V,W)$ by
$$
M(V,W):= (\oplus_{h\in H} \Homc (V_{\out (h)},V_{in(h)}) \oplus
(\oplus_k \Homc(W_k, V_k) \oplus \Homc ( V_k,W_k))
$$. 
For an element of $M(V,W)$ we denote its components by $B_h, i_k,j_k$.\par
Notice  $M(V,W)$ is very close to the vector space $X_V$ defined in
\S2.3. In fact when $W_k=0, \forall k\in I$,  $M(V,0)= X_V$.\par
As in \S2.3, we can define a symplectic form $\omega_{\Bbb C}$ on
 $M(V,W)$ by:
$$
\omega_{\Bbb C} ((B,i,j),(B',i',j')):= \sum_{h\in H}tr(\epsilon(h) B_h
B_{\bar h}') + \sum_k tr(i_ki_k' - i_k'j_k)
$$
where $\epsilon(h)= 1 $ if $h\in \Omega$, 
 $\epsilon(h)=- 1 $ if $h\in \bar\Omega$. The symplectic vector space 
 $M(V,W)$ decomposes into the sum $M(V,W) = M_\Omega \oplus  M_{\bar \Omega}$
or Lagrangian subspaces:
\beqn
M_\Omega:= & (\oplus_{h\in \Omega}  \Homc (V_{\out (h)},V_{\inn(h)}) \oplus
(\oplus_k \Homc(W_k, V_k) \\
M_\Omega:= & (\oplus_{h\in \bar \Omega}  \Homc (V_{\out (h)},V_{\inn(h)}) \oplus
(\oplus_k \Homc( V_k, W_k) 
\endeqn

We can consider $ M_{\bar \Omega}$ as a dual space of $M_\Omega$ via
$\omega_{\Bbb C}$, and $M(V,W)$ as the cotangent bundle of $ M_\Omega $.\par
Since $V$ and $W$ are hermitian vector spaces, $\Homc (V,W)$ has a 
hermitian inner product defined by $(f,g) = tr(fg^*)$, where $g^*$
is the hermitian adjoint. We introduce a natural 
hermitian inner product on  $M(V,W)$ induced from the ones on
 $V$ and $W$ in this way. We introduce a quaternion module structure on 
  $M(V,W)$ by the original complex structure $I$ together with a new
complex structure $J$ by the formula $J(m,m') = (-{m'}^*, m^*)$ for
$(m,m') \in M_\Omega \oplus  M_{\bar \Omega}$.\par
The group $G_V = \prod U(V_k)$ acts on  $M(V,W)$ by 
$$
(B_h,i_k,j_k)\rightarrow (g_{in(h)} B_h g^{-1}_{out(h)}, g_k i_k, j_kg_k^{-1})
$$, preserving the hermitian inner product and the $\Bbb H$-module
structure. Let $\mu$ be the corresponding Hyperk\"{a}hler moment map
(cf.[HKLR]). Denote by $\mu_{\Bbb R}, \mu_{\Bbb C}$ its real
and complex components.  Explicit forms of  $\mu_{\Bbb R}, \mu_{\Bbb C}$
can be found in \S2 of [Na1].
Let $g_V$ be the Lie algebra of  $G_V$ and let $Z_V \subset g_V$
denote the center. Choose an element $\zeta = (\zeta_{\Bbb R}, 0) \in
Z_V \oplus (Z_V\otimes \Bbb C)$ with 
$\zeta_{\Bbb R} = ((i/2)\zeta_1,..., (i/2)\zeta_n)$ satisfies the following
condition:
$$
\zeta_i >0, \forall i
$$. We define a  Hyperk\"{a}hler quotient $M_\zeta(V,W)$ of  $M(V,W)$ by $G_V$
as follows:
$$
M_\zeta(V,W):= \{(B,i,j)\in M(V,W)| \mu(B,i,j) = - \zeta \}/G_V
$$
According to the results of [Na1] ,
$M_\zeta(V,W)$ is a smooth  Hyperk\"{a}hler manifold. \par 

When the
underlying graph is an affine or finite A-D-E graph, $M_\zeta(V,W)$
is a 
certain deformation of singlularities of the moduli space of 
instantons over A-L-E spaces with unitary gauge group 
(cf.Prop.2.15 of  [Na1]). In [Na2], $M_\zeta(V,W)$ is also
refered to as the moduli space of certain torsion free sheaves. \par
We will be interested in certain Lagrangian subvarieties of 
$M_\zeta(V,W)$. Let us define:

$$
\Lambda_{V,W} = \Lambda_{V} \times \oplus_k \Homc(V_k,W_k)
$$
where $\Lambda_{V}$ is the variety defined in definition 2.2.1.
Let us define the set of stable points by:
$$
H^s = \{ m\in \mu_{\Bbb C}^{-1}(0)| G_V^{\Bbb C}.m \cap 
\mu_{\Bbb R}^{-1}(-\zeta_{\Bbb R}) \neq \emptyset \}
$$.
An alternative algebraic characterization of 
$H^s$ is given in Prop.3.5 of [Na1].\par
The Lagrangian subvariety of  $M_\zeta(V,W)$ is defined by:
$$
{\cal L}(V,W):=\Lambda_{V,W} \cap H^s / G_V^{\Bbb C}
$$ (cf.Remark 5.11 of [Na1]).\par
Let Irr${\cal L}(V,W)$ (resp.  Irr$\Lambda_{V,W}$) be the 
set of irreducible components of ${\cal L}(V,W)$ (resp. $\Lambda_{V,W}$).
Since $H^s \rightarrow M_\zeta(V,W)$ is a principal
$G_V^{\Bbb C}$ bundle, Irr${\cal L}(V,W)$
can be identified with 
$$
\{ Y\in  Irr \Lambda_{V,W} | Y \cap H^s \neq \emptyset
\}$$. This can be also identified with
$$
\{ Y\in  Irr \Lambda_{V} | Y \times \oplus_k  \Homc (V_k,W_k)
\cap H^s \neq \emptyset
\}
$$.\par
Recall ${\cal K}_{\Omega}$ is the Grothendieck group
defined in \S1.2. Let $ {\cal K}_{\Omega}^1(W)$ be the
${\Bbb Z} [v,v^{-1}]$ submodule of 
 ${\cal K}_{\Omega}$ generated by complexes $L$ in ${\cal Q}_{\Omega}$
such that 
$$
SS(L) \times \oplus_k  \Homc(V_k,W_k) \cap H^s = \emptyset
$$.
The following lemma is proved in essentially the same way 
as the proof of lemma 11.5 of [Na1]:
\begin{lemma}
$ {\cal K}_\Omega^1(W) \otimes {\Bbb Q}(v)$ is a left ideal in 
${\cal K}_\Omega \otimes {\Bbb Q}(v)$ containing $F_i^{W_i +1}$.
\end{lemma}
Lemma 3.1.1 implies 
$$
\{ L_{b,\omega}(V) | SS(L_{b,\omega}(V)) \times \oplus_k  
\Homc(V_k,W_k) \cap H^s = \emptyset \}
\supset \cup_{i\in I} {^\sigma\text{\bf B}_{V;i;\geq W_i +1}}
$$ by Th.14.3.2 and \S14.4.1 of [Lu].
So:
$$
\{ L_{b,\omega}(V) | SS(L_{b,\omega}(V)) \times \oplus_k
\Homc(V_k,W_k) \cap H^s \neq \emptyset \}
\subset \text{\bf B}_{V}(dimW)
$$ where $\text{\bf B}_{V}(W):= \text{\bf B}(dimW) \cap 
{^\sigma\text{\bf B}_{V}}$
and we have identified $dim W$ with the dominant integral weights of
$U$ as in \S2.2.\par
\begin{lemma}

The map $\psi$ from $ \{ \Lambda_b(V) \in Irr \Lambda_V|
\Lambda_b(V) \times \oplus_k
\Homc(V_k,W_k) \cap H^s \neq \emptyset \}$
to  $ \{  L_{b,\omega}(V) | SS(L_{b,\omega}(V)) \times \oplus_k
\Homc(V_k,W_k) \cap H^s \neq \emptyset \}$ is one-to-one and surjective.
Moreover,
 $ \{  L_{b,\omega}(V) | SS(L_{b,\omega}(V)) \times \oplus_k
\Homc(V_k,W_k) \cap H^s \neq \emptyset \} = \text{\bf B}_{V}(dimW)$.
\end{lemma}
\begin{pf}
By Th.2.3.1, 
$$
SS(L_{b,\omega}(V)) \supset \Lambda_b(V)
$$, hence the map $\psi$ as in Lemma 3.1.2 is well defined and is clearly
one-to-one by definition. To show that  $\psi$ is surjective, 
we just have to show that:
\beq
\sharp & \{ \Lambda_b(V) \in Irr \Lambda_V|
\Lambda_b(V) \times \oplus_k
\Homc(V_k,W_k) \cap H^s \neq \emptyset \} \geq \\ 
\sharp & \{  L_{b,\omega}(V) | SS(L_{b,\omega}(V)) \times \oplus_k
\Homc(V_k,W_k) \cap H^s \neq \emptyset \}.
\endeq
 By the remarks after Lemma 3.1.1, we have:
\beq
\sharp & \{  L_{b,\omega}(V) | SS(L_{b,\omega}(V)) \times \oplus_k
\Homc(V_k,W_k) \cap H^s \neq \emptyset \} \leq \\
 & \sharp \text{\bf B}_{V}(dimW).
\endeq

 On the other hand, let us prove that:
\beq
 \ \ \sharp  \{ \Lambda_b(V) \in Irr \Lambda_V|
\Lambda_b(V) \times \oplus_k
\Homc(V_k,W_k) \cap H^s \neq \emptyset \} = 
\sharp \{Irr {\cal L}(V,W) \}
\endeq
Notice that when the quiver is affine or finite A-D-E, (3.5) is 
(2)  of Prop.10.15 in [Na1] which follows from (3) of 4.16 in [L1]. The
exactly same proof with
possible change of notations 
applies to the general case by using (0) of Th.2.3.1
which is a generalization of (3) of  4.16 in [L1].
 Combine the above with  (d) of Th.33.1.3. and
Th.14.4.11 of [Lu], we have:
$$
\sharp \{Irr {\cal L}(V,W) \} = \sharp \text{\bf B}_{V}(dimW) 
$$. 
It follows that all the $\leq$ or $\geq$ above are actually $=$ and
the lemma is proved. 
\end{pf}
Let us denote by $\Lambda_W$ the simple integrable $U$ module 
with dominant integral weight $dimW$. By Lemma 3.1.1, there is a natural
surjective $U^-$-linear map
$$
\Pi: \Lambda_W \rightarrow {\cal K}_{\Omega} \otimes {\Bbb Q}(v)/
 {\cal K}_{\Omega}^1(W) \otimes {\Bbb Q}(v)
$$.  The dimension of the weight space (with weight $W-V$)
 $\Lambda_W(V)$  is $\sharp \text{\bf B}_{V}(dimW)$,
which is the same  as the corresponding
weight space (with weight $W-V$) 
of $ {\cal K}_{\Omega} \otimes {\Bbb Q}(v)/
 {\cal K}_{\Omega}^1(W) \otimes {\Bbb Q}(v)
$  by Lemma 3.1.2. 
It follow that $\Pi$ is an isomorphism. Let us record this result
in the following theorem:
\begin{thm}
$\Pi$ is an isomorphism. The set Irr${\cal L}(V,W)$ parameterizes
a basis of the weight space $\Lambda_W(V)$.
\end{thm}
In the case of finite A-D-E quivers, Th.3.1.3 is proved in \S11 of
[Na1] by a different method.\par
\subsection{ Quiver varieties with symmetries}
Assume that $(I,H)$ is a finite  graph with a
nontrivial admissible automorphism $a$. 
 We shall fix an orientation
$\Omega$ which is compatible with $a$. Such  an orientation
always exists by \S12.1.1 of [Lu].\par
Recall $u^a$ denotes the universal enveloping algebra of
the corresponding non-symmetric Kac-Moody algebra (cf.\S1.4).
For $W\in {\Bbb N}[I]^a \simeq {\Bbb N}[I^a]$, we denote by 
$L^a(W)$ the simple integrable module of  $u^a$ with
integral dominant weight $dimW$ and $L^a(V,W)$ the subspace 
 of $L^a(W)$ as defined for any simple integrable module at the end of
\S2.1. Here 
$V\in {\Bbb N}[I]^a \simeq {\Bbb N}[I^a]$. \par 
The natural action of $a$ on $M(V,W)$, defined similarly as
the action of  $a$ on $X_V$ in \S1.4, induces an isomorphism,
still denoted by $a: M_\zeta (V,W) \rightarrow M_{a(\zeta)}(a(V),a(W))$,
where $a(\zeta) = (a(\zeta_{\Bbb R}), 0)$ and
$(a(\zeta_{\Bbb R}))_i = (\zeta_{\Bbb R})_{a(i)}$.
Let us fix $\zeta$ as in \S3.1 with  $a(\zeta) =\zeta$.
It is easy to see that only in the case  
$W,V\in {\Bbb N}[I]^a \simeq {\Bbb N}[I^a]$
that $a$ maps $M_\zeta (V,W)$ to itself.\par
\begin{thm}
 Let $W,V\in {\Bbb N}[I]^a \simeq {\Bbb N}[I^a]$. Then:\par
dim$L^a(V,W)= \sharp \{ Y\in Irr{\cal L}(V,W) | a.Y=Y \}$
\end{thm}
\begin{pf}
By Lemma 2.4.1 and Lemma 3.1.2, we have:

\beqn
& \sharp \{ Y\in Irr{\cal L}(V,W)|a.Y=Y \} = \\
& \sharp \{ L_{b,\Omega}(V)|   
  SS( L_{b,\Omega}(V)) \times \oplus_k \Homc (V_k,W_k) \cap H^s \neq \emptyset,
a^*  L_{b,\Omega}(V)
=  L_{b,\Omega}(V) \}\\
&= \sharp \{ L_{b,\Omega}(V) \in \text{\bf B}_V(W)| a^*  L_{b,\Omega}(V)
=  L_{b,\Omega}(V)  \}
\endeqn

By Th.14.4.9 and (d) of Th.13.1.3 of [Lu],
$$
\sharp \{ L_{b,\Omega}(V) \in \text{\bf B}_V(W)| a^*  L_{b,\Omega}(V)
=  L_{b,\Omega}(V)  \} = dim L^a(V,W)
$$. Combining this with the previous identities, the theorem is proved.
\end{pf}
We'd like to describe an application of Th.3.2.1 by constructing  
functions with certain modular properties from the geometry of 
 $M_\zeta(V,W)$.  Suppose the underlying graph is an affine A-D-E
graph with a symmetry $a$. For simplicity we shall restrict to the case
where $(I^a,\cdot)$ is the Cartan datum of non-twisted affine Kac-Moody
algebra.\par
Let $V,W\in {\Bbb N}[I]^a \simeq  {\Bbb N}[I^a]$. 
On the Dykin diagram associated with $(I^a,\cdot)$, assume 
$\tilde i$ corresponds to the label $\alpha_i$ on
Page 54 of [Kac].
We identify $W$
with the dominant integral weight of $g^a$ by identifying $\tilde i\in I^a$
with the  fundamental weight $\Lambda_i$ of  $g^a$ (cf. \S12 of
[Kac]).  Define $V_0$ to be  
$V_r$ for any $r\in \tilde 0$. Let $m_W$ be defined as in formula 12.7.5. on
Page 226 of [Kac] .\par
Let $k(V,W)$ be the complex dimension of $M_\zeta(V,W)$.
Define 
$$
L(W) := \bigoplus_{V,k=k(V,W)} H^k(M_\zeta(V,W), {\Bbb Q})
,$$
 
$$
L(V,W):= H^{k(V,W)}(M_\zeta(V,W), {\Bbb Q})
$$. By Th.10.16
of [Na1], $L(W)$ is the simple integrable module of $u$ with 
integral dominant weight $W$.\par
Let $V,W\in {\Bbb N}[I]^a \simeq  {\Bbb N}[I^a]$,
define a map $L_0': L(V,W)\rightarrow  L(V,W)$ by
$L_0'x= (m_W + V_0) x, \forall x\in  L(V,W).$ \par
If $V$ or $W$ is not an element of ${\Bbb N}[I]^a$, we define
$L_0': L(V,W)\rightarrow  L(V,W)$ to be identity map. \par
The isomorphism $a: M_\zeta (V,W) \rightarrow M_{a(\zeta)}(a(V),a(W))$
induces a natural map :
$$a^*:  H^k(M_\zeta(a(V),a(W)), {\Bbb Q}) 
\rightarrow  H^k(M_\zeta(V,W), {\Bbb Q})
$$.
Let $q$ be an indeterminate. We define a weighted map 
$$
a^*(q):= a^* q^{L_0'}: H^k(M_\zeta(a(V),a(W)), {\Bbb Q})
\rightarrow  H^k(M_\zeta(V,W), {\Bbb Q})
$$. Define 
$$
Ch^a(W):=Tr(a^*(q))
$$ where $a^*(q)$ is an operator from 
 $L(W)$ to itself. It is obvious from the definition that 
$Ch^a(W)=0$ unless $W\in {\Bbb N}[I]^a \simeq {\Bbb N}[I^a]$. Moreover,
the only nonzero contribution to $Ch^a(W)$ comes from the term
$a^*(q):  H^k(M_\zeta(V,W), {\Bbb Q})
\rightarrow  H^k(M_\zeta(V,W), {\Bbb Q})
$ with $W,V\in {\Bbb N}[I]^a$. By Cor.5.5 of [Na1] and Th.3.2.1, 
The trace of the map $a^*(q):  H^k(M_\zeta(V,W), {\Bbb Q})
\rightarrow  H^k(M_\zeta(V,W), {\Bbb Q}) 
$ is equal to $q^{m_W + V_0} $ times dim$L^a(V,W)$. 
  It is easy to check that
$L_0'$ is the same as the operator $L_0^g - \frac{1}{24}c(k)$
in formula 13.8.3 on Page 264 of [Kac]
on $L^a(V,W)$.
Hence $Ch^a(W)= \chi_W (\tau,0,0)$, where 
$\chi_W (h,0,0)$ is defined in formula 13.8.3 on Page 264 of [Kac],
and $q=exp(2\pi i \tau)$.
Let us record this result in the following corollary:
\begin{cor}

 $Ch^a(W) = \chi_W (\tau,0,0)$,
where
$\chi_W (\tau,0,0)$ is defined in formula 13.8.3 on Page 264 of [Kac].
\end{cor}
 By  the above corollary, the functions  $Ch^a(W)$ have certain 
modular properties as in Th.13.8 of  [Kac].\par 
In the case of twisted affine algebras, one can define similar
functions as $Ch^a(W)$ from
the geometry of $M_\zeta(V,W)$ such that it coincides with the normalized
character defined on Page 264 of [Kac] when
$z=\mu=0$. Such functions also have certain
modular properties as shown in Th.13.9 of [Kac]. 
\section{Questions}
The proof of Th.3.2.1. is based on the deep theory of canonical bases of
quantized Kac-Moody algebras. It will be very interesting to proceed 
further to construct a natural analogue of Hecke operators in [Na1] on the
set of constructible functions on ${\cal L}^a(V,W)$ such that these
operators define a representation of $u^a$. Th.3.2.1 suggests that it is
possible to do so but we don't know the answer. \par
In [L1], an explicit description of canonical bases for symmetric affine
 Kac-Moody algebras are given. It will be interesting to use
this result to give a proof of Th.2.3.1 in the affine A-D-E case. This will
avoid the use of Th.2.3.1 in the proof of Cor.3.2.2.\par
As mentioned in the introduction, Cor.3.2.2. suggests that it should be
possible to take into account the actions of
the nontrivial admissible automorphism in the formulation of twisted N=4
SYM theory on certain A-L-E spaces such that the  partition functions
are identified with the functions in Cor.3.2.2. Once this is done, Cor.3.2.2.
together the results of [Kac] can be seen as further evidence in 
favor of S-duality.\par
On Page 93 of [VW] , it is asked if one can find for each rational conformal
field (RCFT) theory a possibly non-compact four-manifold
whose N=4 twisted partition function gives the characters of that 
RCFT. By Cor.3.2.2, it seems that the natural candidate in the case of
Wess-Zumino-Witten model based on non-simply-laced group should be certain
 A-L-E spaces together with a  nontrivial admissible automorphism.
However it is not clear to us what exactly the  four-manifold
should be.\par
S-duality of N=4 SYM theory in 4-dimension 
can be derived by dimensional reduction from the duality between
heterotic string compactified on 4-dimensional torus and type IIA
string  compactified on K3 surface
\footnotemark\footnotetext{ However there is a constraint on the   
possible rank of the gauge group allowed, cf. [LLL] }. 
The work of [Na1] has been
interpreted in this string-string duality context in [HM] and [V].
It will be interesting to see if one can take into account the
action of a non-trivial admissible automorphism in the 
formulation of heterotic/ type IIA duality so at least part of
Cor.3.2.2. follows from such a formulation.\par

Finally, let us note that in the theory of duality between F theory
compactified on Calabi-Yau manifolds and heterotic string compactified on
K3 surfaces,  the Dykin diagram automorphisms corresponding to the 
non-trivial admissible automorphism of finite A-D-E graph
have been used to explain the appearance of non-simply laced
gauge group in  heterotic string (cf. [MV] and [A]). Is it possible
to make connections with the results of this paper by the web of
string dualities?

%%%%%%%%%%%%%%%%%%%%%%%%%%%%%%%%%%%%%%%%%%%%%%%%%%%%%%%%%%%%%%%%%%%%%%%%%%%%%%
 

\begin{thebibliography}{[MFT50]}
\bibitem[BBD]{}A.~A.~Beilinson, J.~Bernstein and P.~Deligne,
{\it Faisceaux perverse}, Ast\'{e}risque 100 (1982).
\bibitem[Kac]{} V.~G.~Kac,
{\it Infinite dimensional Lie algebras}, Cambridge Univ. Press. (1990).
\bibitem[K]{} M.~Kashiwara, Y.~Satio, {\it Geometric Construction
of Crystal Bases}, q-alg/9606009.
\bibitem[Lu]{L;qg}{}G.~Lusztig,
{\it Introduction to quantum groups},
Birkh\"{a}user, (1994)
\bibitem[L1]{}G.~Lusztig,
{\it,Affine Quivers and Canonical bases}, Publ.Math\'{e}matiques, N0.76,(1992),
111-163.
\bibitem[Na1]{}H.~Nakajima,
{\it Instantons on ALE spaces, quiver varieties and Kac-Moody algebras},
Duke Math. Journal (1994) 365-416.   
\bibitem[Na2]{}H.~Nakajima,
{\it Instantons and affine algebras}, alg-geom/9510003.
\bibitem[VW]{} C.~Vafa and E.~Witten,
{\it A Strong coupling test of S-duality}, Nucl.Phys.431 (1994,3-77
\bibitem[V]{} C.~Vafa,
{\it Instantons on D-branes}, hep-th/9512078.
\bibitem[LLL]{} K.~Landsteiner, E.~L\'{o}pez and D.A.~Lowe,
{\it Evidence for S-duality in N=4 Supersymmetric Gauge theory},
hep-th/9606146.
\bibitem[KS]{KS}{}M.~Kashiwara and P.~Schapira,
{\it Sheaves on Manifolds},
Grundlehren Math. Wiss. 292, Springer-Verlag, Berlin (1990).
\bibitem[GL]{}I.~Grojnowski and G.~Lusztig,
{\it A comparison of bases of quantized enveloping algebras},
Contemporary Math. 153 (1993) 11--19.
\bibitem[HM]{}J.A.~Harvey and G.Moore,
{\it On the algebra of BPS states}, hep-th/9609017.
\bibitem[MV]{}M.~Bershadsky, K.Intriligator, S.Kachru, D.R.Morrison,
V.Sadov and C.Vafa,
{\it Geometric Singularities and Enhanced Gauge Symmetries},
hep-th/9605200.
\bibitem[A]{}P. S.~Aspinwall, {\it The SO(32) Heterotic String
on a K3 Surface}, hep-th/9605131.
\bibitem[KP]{} V.~Kac and D.~Peterson, {\it Infinite-dimensional Lie
algebras, theta functions and Modular forms}, Adv.Math. 53 (1984), 125-264.
\bibitem[HKLR]{} N.J.~Hitchin, A.~Karlhede, U.~Lindstr\H{o}m and
M.~Ro\v{c}ek, {\it Hyperk\"{a}hler metrics and supersymmetry,}
Comm.Math.Phys.108 (1987),535-589.
\end{thebibliography}
\end{document}